\documentclass[11pt,a4paper,leqno]{article}

\usepackage{fullpage,amsmath,amsthm,amsfonts,amscd,graphicx,
latexsym,amssymb,eucal,curves}
  \usepackage[all]{xy}
\newcommand{\ZZ}{\mathbb{Z}}

\newcommand{\CC}{\mathbb{C}}

\newcommand{\NN}{\mathbb{N}}

\newcommand{\HH}{\mathbb{H}}

\newcommand{\LL}{\mathbb{L}}
\newcommand{\MM}{\mathbb{M}}
\newcommand{\UU}{\mathbb{U}}
\newcommand{\Schur}{\mathbb{S}}
\newcommand{\HHom}{{\mathbb{H}{\rm om}}}
\newcommand{\WW}{\mathbb{W}}

\newcommand{\QQ}{\mathbb{Q}}
\newcommand{\RR}{\mathbb{R}}

\newcommand{\VV}{\mathbb{V}}
\newcommand{\Ker}{{\rm Ker}}

\newcommand{\End}{{\rm End}}
\newcommand{\EEnd}{{\mathbb{E}{\rm nd}}}

\newcommand{\id}{{\rm id}}
\renewcommand{\Re}{{\rm Re\,}}
\renewcommand{\Im}{{\rm Im\,}}

\newcommand{\Jac}{{\rm Jac} \,}

\newcommand{\Gal}{{\rm Gal}}

\newcommand{\Aff}{{\rm Aff}}

\newcommand{\Sym}{{\rm Sym}}

\newcommand{\mrAV}{A}

\newcommand{\Stab}{{\rm Stab}}
\newcommand{\SL}{{\rm SL}}
\newcommand{\SO}{{\rm SO}}
\newcommand{\diag}{{\rm diag}}

\newcommand{\mrC}{M}

\newcommand{\tr}{{\rm tr}}

\newcommand{\ol}{\overline}

\newcommand{\Hg}{{\rm Hg}}

\newtheorem{Defi}{Definition}[section]

\newtheorem{Rem}[Defi]{Remark}

\newtheorem{Prop}[Defi]{Proposition}
\newtheorem{Lemma}[Defi]{Lemma}
\newtheorem{Cor}[Defi]{Corollary}
\newtheorem{Thm}[Defi]{Theorem}

\begin{document}{\large}
\title{Variations of Hodge structures of a Teichm\"uller curve}
\author{Martin M\"oller}
\maketitle

\begin{quote}
{\footnotesize {\bf Abstract}.
Teichm\"uller curves are geodesic discs in Teichm\"uller space
that project to an algebraic curve in the moduli space $\mrC_g$. 
We show that for all $g \geq 2$ Teichm\"uller curves map to the locus
of real multiplication in the moduli space of abelian varieties.
Remark that McMullen has shown that precisely for $g=2$ the locus
of real multiplication is stable under the $\SL_2(\RR)$-action
on the tautological bundle $\Omega\mrC_g$. 
\newline
We also show that Teichm\"uller curves are defined over number
fields and we provide a completely algebraic description of
Teichm\"uller curves in terms of Higgs bundles. As a consequence we show
the absolute Galois group acts on the set of Teichm\"uller curves.
}
\end{quote}

\section*{Introduction}
The bundle $\Omega T_g^*$ over the 
Teichm\"uller space, whose points parametrize pairs $(X,\omega)$
of a Riemann surface $X$ of genus $g$ with a Teichm\"uller marking
and a non-zero holomorphic $1$-forms $\omega$  has a natural 
$\SL_2(\RR)$-action. Its can be described by postcomposing the charts 
given by integration of $\omega$  with the linear
map and providing the Riemann surface with the new complex
structure. The linear action on the real and imaginary part of
$\omega$ yields a one-form that is holomorphic for the new
complex structure.
\newline
The orbit of $(X,\omega) \in \Omega T_g^*$ projected
to $T_g$ is a holomorphic embedding of the upper half plane into
Teichm\"uller space
$$ \tilde{\j}: \HH \to T_g,$$
which is totally geodesic for the Teichm\"uller metric.
Only rarely these geodesics project to algebraic curves 
$C = \HH/\Stab(\tilde{\j})$ in moduli space $\mrC_g$. These curves are 
called {\em Teichm\"uller curves} and $(X,\omega)$ a
{\em Veech surface}. Coverings of the torus ramified over only
one point give Veech surfaces. Veech (\cite{Ve89})
constructed a series of examples in all genera that do not
arise via coverings. We will say
more about $\Stab(\tilde{\j})$ in Section $1$.
\par
McMullen has shown in \cite{McM1} that in genus $2$ a Teichm\"uller
curve always maps to the locus of real multiplication in the
moduli space of abelian surfaces $\mrAV_2$ or to the locus of
abelian surfaces that split up to isogeny. 
Moreover he has shown that the locus of eigenforms of real 
multiplication is $\SL_2(\RR)$-invariant. Since the
property of having a double zero is also $\SL_2(\RR)$-invariant
he constructs infinitely many examples in $g=2$ which are 
not obtained via coverings of elliptic curves. 
\newline
McMullen also shows that for $g>2$ the locus of
eigenforms is no longer $\SL_2(\RR)$-invariant.
\par 
We investigate the variation of Hodge structures (VHS) of the
family of Jacobians over a Teichm\"uller
curve $C$ or more precisely over an unramified cover of $C$
where the universal family exists (see Section \ref{univfam}).
\newline 
It turns out that the underlying local system has
a $\QQ$-factor which is the sum of $r$ rank $2$ local systems
given by the action of $\Stab(\tilde{\j})$ and its Galois conjugates.
Here $r \leq g $ equals the degree of the trace field of $\Stab(\tilde{\j})$ 
over $\QQ$. As a consequence we obtain for all genera:
\par
{\bf Theorem \ref{Endallgemein}.} {\rm 
The image of a Teichm\"uller curve $C \to A_g$ is contained in the
locus of abelian varieties that split up to isogeny
into $A_1 \times A_2$, where $A_1$ has dimension $r$ and
real multiplication by the trace field $K= \QQ(\tr(\Stab(\tilde{\j}))$,
where $r=[K:\QQ]$. 
The generating differential $\omega \in H^0(X,\Omega^1_X)$
is an eigenform for the multiplication by $K$, i.e.\
$K\cdot \omega \subset \CC \omega$.}
\par
For further results stated in Sections $4$ and $5$ 
it is convenient also to have the language of Higgs bundles
availabe. The basic notions will be recalled in Section $3$.
Simpson's correspondence allows us to go back and forth
from VHS to Higgs bundles. 
\newline
Our study of Teichm\"uller curves is inspired by the work
of Viehweg and Zuo on Shimura curves (\cite{ViZu04}). They
show that the VHS over a Shimura curve is build out of sub-VHS, 
whose Higgs field is 'maximal' and unitary sub-VHS. 
Maximality here means that the Higgs field is an isomorphism.
The notion is motivated since this happens if and only
if a numerical upper bound (Arakelov inequality) is attained.
\par
Teichm\"uller curves are in some sense similar to Shimura
curves, since their VHS always contains a sub-VHS, that
is maximal Higgs, see Prop.\ \ref{splitasLvar} and its
converse Theorem \ref{conversekrit}. 
\newline
In fact we characterize Teichm\"uller curves algebraically 
using this notion:
\par
{\bf Theorem \ref{algchar}.} {\rm
Suppose that the Higgs bundle of family of curves 
$f: {\cal X} \to C = \HH/\Gamma$ has a rank two 
Higgs-subbundle with maximal Higgs field. 
\newline
Then $C \to \mrC_g$ is a finite covering of a Teichm\"uller curve.
}
\par
Suppose that we have a Teichm\"uller curve with $r=g$. This
implies that it does not arise from lower genus via a
covering construction. The above algebraic description implies
that Teichm\"uller curves are defined over number fields
and that $G_\QQ$ acts on the set of Teichm\"uller curves with $r=g$
(see Thm.~\ref{nodefo} and Cor.~\ref{GQacts}).
\par
There are also properties that are significantly different
between Teichm\"uller- and Shimura curves. For example 
consider the locus of real multiplication by some field
in the moduli space of abelian threefolds. There are
Shimura subvarieties properly contained in this locus
that contain (lots of) Shimura curves. But for 
Teichm\"uller curves we show: 
\par
{\bf Theorem \ref{Shimura}.} {\rm
Suppose that $C \to M_g$ is a Teichm\"uller curve with 
$r = g$. 
\newline
Then the Shimura variety parametrizing abelian varieties with
real multiplication by $K$ is the smallest Shimura subvariety
of $\mrAV_g$ the Teichm\"uller curve $C$ maps to.
}
\par 
The above restrictions to the case $r=g$ may be explained as
follows: The VHS over a Teichm\"uller curve has a rank $2$
sub-VHS and its Galois conjugates, which we have control of. 
For the 'rest', which may be enlarged by covering constructions,
we can hardly say anything. It would be intereresting to 
characterize the local systems that can appear in that 'rest'
when $r<g$.
\par
The author thanks Eckart Viehweg for his constant encouragement
to describe the VHS of Teichm\"uller curves and for many
fruitful discussions. He also thanks the referee for
several suggestions.
\par
\section{Setup}

\subsection{Flat structures}
Start with a Riemann surface $X$ of genus $g \geq 1$ and a non-zero 
holomorphic quadratic differential $q \in H^0(X,\Omega_X^{\otimes 2})$. 
This pair determines a {\em flat structure} on $X$: Cover $X$ minus
the zeros of $q$ by simply connected open sets $U_i$ on which a 
$q$ admits a square root $\omega$. Integration of $\omega$ gives
charts $\varphi_i: U_i \to \CC$ such that the transition functions 
$\varphi_j \circ \varphi_i^{-1}$ are translations composed with $\pm 1$.
\newline
Although the $SL_2(\RR)$-action and Teichm\"uller geodesics (see below)
are defined for all pairs $(X,q)$, we consider here only the case where 
$q=\omega^2$ is globally the square of a holomorphic one-form 
$\omega \in H^0(X,\Omega_X)$. This restriction is not too serious since 
a pair $(X,q)$ with $q$ non-square admits a double covering 
$\pi: X' \to X$ branched over the zeros of $q$ of odd order such that
$\pi^* q$ is a square. On the other hand it seems necessary, since
quadratic differentials are not directly reflected in the Hodge structure,
see Section \ref{Splitting}.
\par

\subsection{$\SL_2(\RR)$-action}
Denote by $\Omega T_g$ the vector bundle over the
Teich\-m\"uller space $T_g$ whose fibres are global holomorphic
$1$-forms. Equivalently, $\Omega T_g$ is the total space of the
pullback of the tautological bundle (hodge bundle) over $M_g$
to $T_g$. Let $\Omega T_g^*$ denote the bundle with the zero section
removed. $\Omega T_g^*$ carries a natural $\SL_2(\RR)$-action: 
Given $(X,\omega) \in \Omega T_g^*$
and $ A= { \left( \begin{array}{ll} a & b \\ c & d \end{array} 
\right)\ } \in 
\SL_2(\RR)$ consider the $1$-form
\begin{equation}\label{SLaction}
 \eta = \left( \begin{array}{ll} 1 & i \end{array}\right) 
\left( \begin{array}{ll} a & b \\ c & d \end{array} \right)
\left( \begin{array}{l} \Re \omega \\ \Im \omega  \end{array} \right). 
\end{equation}
There is a unique complex structure on the topological surface underlying
$X$ such that $\eta$ is holomorphic. Call this surface $Y$ and
let $ A \cdot (X,\omega) = (Y, \eta)$. Equivalently we may postcompose
the charts $\varphi_i$ of the flat structure by the action of $A$
on $\CC \cong \RR^2$. The new complex structures glue, since 
the transition functions are translations.
\newline
The fibres of the projection
$\Omega T_g^* \to T_g$ are stabilized by $\SO_2(\RR)$.
Consider the quotient mapping $\tilde{\j}: \HH \cong \SL_2(\RR) / \SO_2(\RR) \to  T_g$.
This map is holomorphic and a geodesic embedding for the
Teichm\"uller (equivalently: Kobayashi) metric on $T_g$. In fact
we can identify $\HH$ with the unit disk $t \ol{\omega}/\omega$
for $|t|<1$ in the space of Beltrami differentials, see
\cite{McM1} Section $2$.
\par

\subsection{Teichm\"uller curves}
In some rare occasions the projection $\HH \to T_g \to \mrC_g$ 
to the moduli space of curves of genus $g$ is closed for the
complex topolgy. By \cite{SmWe04} Prop.~8 this is 
equivalent to the setwise stabilizer 
$\Stab(\tilde{\j}):= \Stab_{\Gamma_g}(\tilde{\j})$ of 
$\HH$ under the action of the mapping class group $\Gamma_g$ 
being a lattice in 
$\SL_2(\RR)$. In this case we call $\tilde{\j}$ 
(or $C:=\HH/\Stab(\tilde{\j})$ or $j: C \to \mrC_g$) a {\em Teichm\"uller
curve}. We also say that $(X,\omega)$
{\em generates} a Teichm\"uller curve. If we change the basepoint i.e.\ 
if we start with  $A \cdot (X,\omega)$ instead of
$(X,\omega)$, this amounts to replacing $\Stab(\tilde{\j})$ by
$A^{-1} \cdot \Stab(\tilde{\j}) \cdot A$. Note that $C$ is 
a complex algebraic curve, but it is never complete by \cite{Ve89} 
Thm.~2.1.
\par
\subsection{Universal family} \label{univfam}

We want to consider a Teichm\"uller curve $C \to M_g$ as a
family of curves parametrized by $C$. All the results in the sequel will 
not depend on taking a subgroup $\Gamma$ of finite index 
in $\Stab(\tilde{\j})$.
There is subgroup $\widetilde{\Gamma_g}$ of finite index 
in the mapping class group which is torsion free. 
Hence there is a universal family of
curves over $T_g/\widetilde{\Gamma_g}$. We may take for this purpose 
$T_g/\widetilde{\Gamma_g} = M_g^{[3]}$ the moduli space of curves with
level-$n$-structure for $n = 3$. We lift the Teichm\"uller curve to
this covering: Let $\Gamma_1:= \Stab_{\widetilde{\Gamma_g}}(\tilde{\j})$ 
and let $C^{[3]} = \HH/\Gamma_1 \to  M_g^{[3]} $, be the induced quotient of 
$\tilde{\j}$. We let $f^{[3]}: {\cal X}^{[3]}_{\rm univ} \to M_g^{[3]}$ be
this universal family and we pull it back to $C^{[3]}$.
\newline
This is almost the family of curves we want. Its monodromies around the
cusps are quasi-unipotent (see e.g.\ \cite{Sd73} Thm.~6.1). 
For technical purposes we pass to a finite index 
subgroup $\Gamma$ of $\Gamma_1$ such that there the family of
curves admits a stable model over $\ol{\HH/\Gamma}$ and hence
monodromies around the cusps are unipotent. 
To ease notation we denote the pullback of the universal family 
to this covering by $f: {\cal X} \to C$ and call it a family
of curves {\em coming from a Teichm\"uller curve}. We hope that no confusion
with the original definition of $C$ lying inside $M_g$ will arise.
\newline
The whole situation is summarized in the following diagram: 
$$ \xymatrix{
\ol{{\cal X}} \ar[d]^f   & {\cal X}\ar[d]^f \ar@{_(->}[l] \ar[r]  
& {\cal X}_C^{[3]} \ar[r] \ar[dd]^{f^{[3]}} 
& {\cal X}^{[3]}_{\rm Univ} \ar[dd]^{f^{[3]}} \\
\ol{C} & C:= {\HH/\Gamma} \ar@{_(->}[l] \ar[dr] & & \\
&  & C^{[3]}=\HH/\Gamma_1 \ar[r]\ar[d]  & M_g^{{[3]}} \ar[d] \\
& & C = \HH/\Stab(\tilde{\j})\ar[r] & M_g \\
}$$
Here $\ol{C}$ is a smooth completion of $C$ and $\ol{{\cal X}}$ is a family of
stable curves over $\ol{C}$. 
\newline
We will denote throughout by $g: \Jac({\cal X}/C) \to C$
the family of Jacobians of $f: {\cal X} \to C$ and abbreviate
$A:= \Jac({\cal X}/C)$. 
\par

\subsection{The affine group}
Let $\Aff^+(X,\omega)$ denote the group of orientation-preserving
diffeomorphisms, which are affine with respect to the flat
structure defined by $\omega$. The 'derivative' 
$D: \Aff^+(X,\omega) \to \SL_2(\RR)$ associates 
with such a diffeomorphism its matrix part. We denote
the image of $D$ by $\SL(X,\omega)$. This group is related
(see \cite{McM1} Prop.\ 3.2) to the stabilizer by
\begin{equation} \label{affstab}
\Stab(\tilde{\j}) = \diag(1,-1) \cdot \SL(X,\omega) 
\cdot \diag(1,-1). 
\end{equation}   
The kernel of $D$ is the group of conformal automorphisms of $X$
preserving $\omega$ and corresponds to the pointwise stabilizer
of $\HH$ in $T_g$. 
\par

\subsection{Real multiplication} \label{RM}
Let $(\pi:A_1 \to S, \lambda: A_1 \to A_1^{\vee})$ be family of 
polarized abelian varieties of dimension $r$ over some base $S$
over the complex numbers.
This family is said to have {\em real multiplication by a totally real number 
field $K$} with $[K:\QQ]=r$,
if there is an inclusion of $\QQ$-algebras
$$ K \hookrightarrow \End(A_1/S) \otimes \QQ.$$
Since we are working in characteristic zero, we do not need
to impose supplementary conditions on the tangent space 
in order to have a well-behaved moduli functor, see \cite{vG87}
Lemma X.1.2. Two families of abelian varieties with real 
multiplication are said to be isomorphic, 
if there is an isomorphism respecting the polarization and the $K$-action.
\par
Suppose $S$ is a point. The action of $K$ on $H^0(A_1,\Omega^1_{A_1})$ 
is diagonalizable and we may choose an eigenbasis 
$\{\omega_1,\ldots,\omega_r\}$.
We call $\omega_i$ an eigenform for real multiplication.
Furthermore we may choose a symplectic basis $\{a_i,b_i\} \in H_1(A_1,\RR)$ 
adapted to $\{\omega_i\}$. By definition adapted means that 
$\int_{a_i} \omega_j = 0$ and
$\int_{b_i} \omega_j = 0$ for $i \neq j$ and that 
$ \int_{a_i} \omega_i /\int_{b_i} \omega_i \in \HH$.
\newline
If $S$ is simply connected this holds in the relative setting:
we find sections $\omega_i$ of $\pi_* \Omega^1_{A_1}$ of eigenforms
and sections $a_i, b_i$ of the local system $R^1 \pi_* \RR_{A_1}$ 
that are fibrewise adapted in the above sense.
\newline
This shows that the moduli functor (in the complex category)
for abelian varieties with real multiplication by 
$K$ with an adapted symplectic basis is represented by $\HH^r$.
\newline
If we fix a level structure instead of a symplectic basis and if
we fix the precise type of the real multiplication, i.e.\ an embedding
of an order ${\mathfrak o} \subset K$ into the endomorphism
ring of the abelian variety, two elements of $\HH^r$ become
identified, if they lie in the same orbit of a group $\Delta$
commensurable with $\SL_2({\mathfrak o})$ depending on the level
structure chosen.
Here the $r$ embeddings $K \hookrightarrow \RR$ induce a 
map $\SL_2({\mathfrak o}) \hookrightarrow \SL_2(\RR)^r$ 
and $\SL_2({\mathfrak o})$ acts on $\HH^r$ via this map. 
As above we do not care about specifying the level structure or
the order $\mathfrak o$ since our results will not depend on passing 
to unramified covers.
\newline
For more details see \cite{LaBi92}, \cite{vG87} and \cite{McM1}, \cite{McM3}.
\par
\section{Splitting of the local system} \label{Splitting}

We start with a review on variations of Hodge structures
(VHS) for a family $f: {\cal X} \to C$ of curves over 
a base curve $C$ and we suppose that it extends to
a family of stable curves $\ol{f}: \ol{{\cal X}} \to \ol{C}$. Fix a base point 
$c \in C(\CC)$. The cohomology $H^1(X,\ZZ)$ of the fibre $X$ of $f$ over $c$
is acted upon by $\pi_1(C,c)$. This data is equivalent to having a $\ZZ$-local
system $R^1 f_* \ZZ_{\cal X}$ on $C$, where $\ZZ_{\cal X}$ is the
constant sheaf $\ZZ$ on ${\cal X}$. The associated rank $2g$ vector bundle
$R^1 f_* \ZZ_{\cal X} \otimes_\ZZ {\cal O}_C $ has a distinguished extension
to $\ol{C}$ due to Deligne (see \cite{De70}), 
which we denote by $E:= (R^1 f_* \ZZ_{\cal X})_{{\rm ext}}$. The holomorphic
one-forms form a subspace of $H^1(X,\CC)$. In the relative situation this 
fits together to the bundle $f_* \omega_{\ol{\cal X}/\ol{C}} \subset E$. 
Hodge theory plus Serre duality implies that the quotient
$$E/ f_* \omega_{\ol{\cal X}/\ol{C}} \cong 
(f_* \omega_{\ol{\cal X}/\ol{C}})^\vee \cong
R^1 f_* {\cal O}_{\cal X}.$$
\par
Furthermore the fibrewise symplectic pairing $(\cdot,\cdot)$ on 
$H^1(X,\RR)$ is gives a hermitian form
$$ H(v,w) := i (v,\ol{w}). $$ 
The data $(R^1 f_* \ZZ_{\cal X}, \; 
f_* \omega_{\ol{\cal X}/\ol{C}} \subset E, \; H)$ 
form a polarized VHS (pVHS) of weight $1$. We are also interested
in sub-VHS. For this purpose we have to allow the underlying 
local system to be only a $F$-local system for some field 
$F\subset \CC$. In this case we talk of an $F$-VHS.
\par
Weight one will be sufficient for a reader mainly intereted in understanding
the results on Teichm\"uller curves in this section. In general and for
Sections \ref{Higgssection} and \ref{noconstant} we recall that a polarized  
$L$-VHS of weight $n$ over a curve $C$ consists of 
\newline
$\bullet$ a $L$-local system $\VV_L$ over $C$ for a field $L\subset \CC$
which gives
rise to a connection $\nabla$ on $E := 
(\VV_L \otimes_L {\cal O}_C)_{{\rm ext}}$ and
\newline
$\bullet$ a decomposition $E=\oplus_{p\in\ZZ} E^{p,n-p}$ into
$C^{\infty}$-bundles, such that 
\begin{itemize}
\item[i)] $E^p := \oplus_{i\geq p} E^{i,n-i}$ are holomorphic subbundles
and $\ol{E^q} := \oplus_{i\leq n-q} E^{i,n-i}$ are antiholomorphic
subbundles and
\item[ii)] $\nabla(E^p) \subset \Omega^1_C \otimes E^{p-1}$ and 
$\nabla(\ol{E^q}) \subset \Omega^1_C \otimes \ol{E^{q-1}}$.
\end{itemize}
The VHS is {\em polarized} if there exists a hermitian form
$H(\cdot,\cdot)$ on $E$ such that the decomposition $E=\oplus E^{p,n-p}$
is orthogonal and such that $(-1)^{p-q}H(v,\ol{v})>0$ for all $v \in E^{p,q}$.
\newline
A homomorphisms of VHS $\varphi:E \to F$ of bigree $(r,s)$ is a 
homomorphism of the underlying local systems such that
$\varphi(E^{p,q}) \subset F^{p+r, q+s}$.
\par
We now specialize to the case of a family $f$ coming from a Teichm\"uller
 curve as in Section \ref{univfam}:
\newline 
Consider the cohomology of the fibre $X$ over $c \in C(\CC)$ with $\RR$
coefficients. Since $C = \HH/\Gamma$, the subspace $\LL_{c} := 
\langle {\Re} \omega, {\rm Im} \omega \rangle_\RR \subset H^1(X,\RR)$ 
is invariant under the action of
$$\pi_1(C,c) = \Gamma \subset \SL_2(\RR)$$ 
by the defining equation (\ref{SLaction}). 
\newline
We thus obtain a rank $2$ linear subsystem
$$\LL_\RR \subset \VV_\RR = (R^1 \pi_* \QQ_{\cal X}) \otimes_\QQ \RR, $$
\par
whose fibre over $c$ is $\LL_{c}$. We now want to decompose the
VHS and we recall Deligne's semisimplicity for convenience. 
For our purpose the more well-known result
for $\QQ$-VHS is not sufficient, we need it for $\CC$-VHS.
\par
\begin{Thm} \label{Del87prop} (Deligne, \cite{De87} Prop.\ 1.13)
The local system $\VV_\CC$ on the algebraic curve $C$ decomposes as 
$$\VV_\CC = \bigoplus_{i=1}^n (\LL_i \otimes W_i),$$
where $\LL_i$ are pairwise non-ismorphic 
irreducible $\CC$-local systems and $W_i$
are non-zero $\CC$-vector spaces.
\newline
Moreover the $\LL_i$ and the  $W_i$ carry 
polarized VHS, whose tensor product and sum gives back the Hodge structure 
on $\VV_\CC$. The Hodge structure on the $\LL_i$ (and $W_i$ 
is unique up to a shift of the bigrading.
\end{Thm}
\par
Deligne states the above proposition for weight zero, but the
proof carries over to any other weight without changes.
\newline
A few words on the proof: The semisimplicity of $\VV_\CC$
is due to the existence of a polarisation, see \cite{De71} Thm.~4.2.6.
Given a local system $\VV$ on $C$ we denoote by 
$\End(\VV_\CC) = H^0(C, \EEnd(\VV_\CC))$ the global sections
of the local system of endomorphisms. Then, due to
the irreducibility of the $\LL_i$, we have
$$\End(\VV_\CC) = \prod_{i=1}^n \End(W_i).$$
Now the components of the Hodge decomposition of $\End(\VV_\CC)$
are again flat global sections. This is the point where one
needs that $C$ is an algebraic curve or at least that any 
subharmonic function on $C$ bounded above is constant. 
We thus obtain a Hodge structure on $\prod_{i=1}^n \End(W_i)$.
The remaining steps consist of showing that this Hodge structure
comes from one on $W_i$ and to use suitable projections $\VV_\CC \to \LL_i$
to provide the $\LL_i$ with a VHS.
\par 
One more notation: Let $L$ be a field. An $L$-local system 
$\WW$ on $C$ is {\em defined over} a subring $R \subset L$  
if there is a local system of torsion-free $R$-modules $\WW_R$ 
with $\WW \cong \WW_R \otimes_R L$. Equivalently the representation
of $\pi_1(C,c)$ of the fibre $\WW_c \cong L^{{\rm rank}\, \WW}$ 
is conjugate to a representation that factors though
${\rm GL}(R^{{\rm rank}\, \WW})$.
\par
We come back to the case of Teichm\"uller curves:
\newline
Since the Fuchsian group $\Gamma$ defines a local
system that is irreducible over $\CC$, we may suppose that
$\LL_1 = \LL_\CC$. Let $K = \QQ(\tr(\Gamma))$ denote the
trace field of the Fuchsian group. It is a real number field
of degree $r:=[K:\QQ] \leq g$ over $\QQ$ (\cite{McM1} Thm..1).
\par
\begin{Lemma} \label{fromTake}
$\LL_\RR$ is defined over number field $K_1 \subset \RR$,
which has degree at most two over $K$.
\end{Lemma}
{\bf Proof:} By Takeuchi (\cite{Ta69}) it suffices to take
$K_1 = K(\lambda)$, where $\lambda$ is an eigenvalue of a 
hyperbolic element in $\Gamma$. \hfill $\Box$
\par
We denote by $\LL$ the $K_1$-local system such that
$\LL \otimes_{K_1} \RR = \LL_\RR$. Denote by $L$ the 
Galois closure of $K_1/\QQ$ and let $r=[K:\QQ]$. 
The $\Gal(L/\QQ)$-conjugates $\LL^\sigma$ of $\LL$ are also 
rank $2$ irreducible local subsystems of $\VV_L$. Hence the
$(\LL^\sigma)_\CC$ appear among the $\LL_i$ in the 
semisimplicity theorem. 
\par
\begin{Lemma} \label{W1dim}
For $\sigma \in \Gal(L/\QQ)$ the local
subsystems $\LL$ and $\LL^\sigma$ in $\VV_L$ are isomorphic if and only if
$\sigma$ fixes $K$. Moreover $W_1$ is one-dimensional. 
\end{Lemma}
\par
{\bf Proof:}  If $\LL$ and $\LL^\sigma$ are isomorphic,
the traces of $\Gamma$ are invariant under $\sigma$, hence $\sigma$ 
has to fix $K$. For the converse 
let $\phi \in \Aff(X,\omega)$ be a diffeomorphism, whose
image $D(\phi)$ lies in the subgroup $\Gamma$ of $\SL(X,\omega)$
and which is hyperbolic, i.e.\ whose trace $t$ has absolute value 
greater than $2$. We let $\phi^*$ denoted the induced morphism on 
$H^1(X,\QQ)$. By \cite{McM2} Thm.~9.5 the sub-vectorspace
$\Ker((\phi^* + (\phi^*)^{-1}) - t\cdot\id)$ of $H^1(X,\RR)$
is two-dimensional and contains $\LL_{c}$. If 
$\LL \not\cong \LL^\sigma$ for $\sigma \in \Gal(L/K)$ or if
$\dim W_1 \geq 2$ then this subspace would have to be of
dimension greater than two.
\hfill $\Box$
\par
We denote the conjugate $L$-local systems of $\LL$ under a fixed 
system of representatives of $\Gal(L/\QQ)/\Gal(L/K)$ by 
$\LL = \LL_1,\ldots,\LL_r$ with the slight abuse of notation that
they appear as $\CC$-local systes in the semisimplicity theorem. 
\newline
We can now sum up: 
\par 
\begin{Prop} \label{splitasLvar}
The local system $\VV = R^1 f_*\ZZ$ splits over
$\QQ$ as
$$ \VV_\QQ = \WW_\QQ \oplus \MM_\QQ, \quad {\rm where} \quad 
\WW_L = \LL_1 \oplus \cdots \oplus \LL_r, $$
such that each of the $\LL_i$ carries a polarized $L$-VHS
of weight one and $\MM_\QQ$ carries a polarized $\QQ$-VHS  
of weight one whose sum is the VHS on $\VV$. Moreover none
of the $(\LL_i)_\CC$ is contained in $\MM_\CC$.
\end{Prop}
\par
{\bf Proof:} Apply Thm. \ref{Del87prop} and let 
$$\MM_\CC := \sum_{i=r+1}^n (\LL_i \otimes W_i).$$ 
The vector spaces $W_i$ for $i=2,\ldots,r$ are one-dimensional, 
since otherwise by Galois conjugation we would have 
$\dim W_1 \geq 2$ contradicting Lemma \ref{W1dim}.
The local system $\WW_L$ is defined over $\QQ$, since $\Gal(L/\QQ)$
acts by permuting its summands. Hence also the complement $\MM_\CC$
is defined over $\QQ$.
\hfill $\Box$.
\par
\begin{Rem} \label{McM42Rem} {\rm 
The argument, which implies that $\dim W_1 = 1$ in the above Lemma
\ref{W1dim} may be explained as follows: The hyperbolic element
$\phi$ is a pseudo-Anosov diffeomorphism of $X$. It has unique
contracting and expanding foliatons and unique maximal and
minimal eigenvalues when acting of $H^1(X,\RR)$. The subspace
singled out by $\LL_1$ is the sum of the corresponding eigenspaces,
compare \cite{McM1} Thm.~5.3. This lemma is the only point,
besides the definition of a Teichm\"uller curve that gives rise
to $\LL_1$, where actually Teichm\"uller theory is used. A 
corresponding algebraic argument seems not available. It would permit
e.g.\ to extend some of the results in \cite{ViZu05} to higher
dimension.
\newline
There is another argument proving $\dim W_1=1$ also using
the Teichm\"uller metric, based on \cite{McM1} Thm.~4.2: 
\newline
Suppose the dimension was $s \geq 2$. Choose sections
$\omega_i(\tau)$ for $i=1,\ldots,s$ of the pullback of
the $(1,0)$-part of  $(\LL_1 \otimes W_1)\otimes_\CC {\cal O}_C$ to 
the universal covering $\HH$. We may suppose that $\omega_1$
generates the Teichm\"uller curve.
Choose a symplectic basis $\{a_i,b_i\}$ for
$i=1,\ldots,s$ of sections of $R_1 f_* \ZZ_{{\cal X}}$ such
that fibrewise $\int_{a_i} \omega_j = \delta_{ij}$ and 
$\int_{b_i} \omega_j = 0$ for $i \neq j$. 
\newline
By construction the map 
\begin{equation} \label{classmap}
\HH \ni \tau \mapsto \int_{b_1} \omega_1(\tau) \in \HH
\end{equation} 
is an isometry. Since the local system $\LL_1$ appears several
the same is true if we replace $\omega_1$ by $\omega_j$ for
$j=2,\ldots,s$ in equation (\ref{classmap}). This contradicts 
the theorem just cited, which states
that these maps are strict contractions for differentials 
not proportial to $\omega_1$.
}\end{Rem}
\par
We now study the implications of the decomposition in 
Prop.\ \ref{splitasLvar} on the endomorphism ring of the family of
Jacobians $A/C$. 
Given a local system $\VV$ on $C$ we denote by 
$\End(\VV) = H^0(C, \EEnd(\VV))$ the global sections
of the local system of endomorphisms of $\VV$. We have 
$$ \EEnd(\VV_L) \supset \EEnd(\WW_L) = \bigoplus_{i=1}^r \EEnd(\LL_i) \oplus 
\bigoplus_{1 \leq i,j \leq r, i\neq j} \HHom(\LL_i,\LL_j). $$
For $a \in K$ consider the following
element of $\End(\VV_L)$: We define its action
on $\LL_i$ as $a^{\sigma_i}\cdot \id$  (where $\sigma_i \in \Gal(L/\QQ)$
maps $\LL_1$ to $\LL_i$). This endomorphism is $\Gal(L/\QQ)$-invariant
hence in $\End(\VV_\QQ)$, and of bidegree $(0,0)$. By Deligne's description of 
endomorphisms of abelian varieties (\cite{De71} Rem.~4.4.6) this
means that the field $K$ is contained in $\End(A) \otimes_\ZZ \QQ$. 
The decomposition of $\VV_\QQ$ over $\QQ$ translates into isogeny of abelian
schemes over $A/C$ 
$$ A \hookrightarrow A_1 \times_C A_2,$$
where $\dim A_1 =r$. 
\par 
By the classification of endomorphism rings of abelian varieties
(see e.g.\ \cite{De71} 4.4.5) a field $K \subset 
\End(A) \otimes_\ZZ \QQ$ is totally real or an imaginary quadratic
extension of a totally real field. Since $K$ arose as a trace field, it
is real and we have shown:
\par
\begin{Prop}\label{LinRR}
The trace field $K=\QQ(\tr(\Gamma))=\QQ(\tr(\SL(X,\omega))$ is totally real.
\end{Prop}
\par
\begin{Thm} \label{Endallgemein}
The image of a Teichm\"uller curve $C \to A_g$ is contained in the
locus of abelian varieties that split up to isogeny
into $A_1 \times_\CC A_2$, where $A_1$ has dimension $r$ and
real multiplication by the field $K= \QQ(\tr(\SL(X,\omega)))$. 
The generating differential $\omega \in H^0(X,\Omega^1_X)$
is an eigenform for multiplication by $K$, i.e.\
$K\cdot \omega \subset \CC \omega$.
\end{Thm}
\par
\begin{Rem}{\rm 
McMullen noted in \cite{McM1} Thm.~7.5 that the eigenlocus for
real multiplication is for $g \geq 3$ no longer invariant
under the action of $\SL_2(\RR)$ on holomorphic differentials
of Riemann surfaces of genus $g$. He starts with an
eigendifferential $\omega = xdx/y$ for $\xi = \zeta_7 + \zeta_7^{-1}$ 
on the curve $X: y^2 = x^7-1$ and shows the following: 
The automorphisms $(x,y) \mapsto (\zeta_7 x,y)$ of 
$X$ are elements of $\Aff^+(X,\omega)$. The 
group $\SL(X,\omega)$ contains an element $\gamma$ with 
trace $\xi$, but nevertheless the $\SL_2(\RR)$-orbit of
$(X,\omega)$ leaves the locus of real multiplication.
This is no contradiction. The point is that the proof of 
Thm.~\ref{Endallgemein} relies on Deligne's semisimplicity, 
which holds only if on the base of the family any subharmonic function
bounded above is constant. The upper half plane 
(also mod $\langle \gamma \rangle$)
is not of this type. 
}\end{Rem}
\par
\begin{Rem}{\rm On the curve $X_n: y^2=x^n-1$ a basis of
holomorphic differentials is $\omega_i=x^{i-1}dx/y$ for $i=1,\ldots,g$.
In \cite{Ve89} Veech showed that the geodesics generated
by $(\omega_1)^2$ are Teichm\"uller curves. Thm.~\ref{Endallgemein}
together with McMullen's theorem shows that, say for odd $n$,
the geodesics generated by $(\omega_i)^2$ for $1<i<g$ 
are not Teichm\"uller curves. In fact one can write
$(\omega_i)^2 = \omega_{i-1} \omega_{i+1}$ and apply 
Ahlfors variational formula to show that the $\SL_2(\RR)$-deformations
leave the locus of real multiplication.
}
\end{Rem}
\par

\subsection{The quaternion algebra assoicated with $\SL(X,\omega)$}
For a Fuchsian group $\Gamma$ contained in $\SL_2(K_1)$ for
some number field $K_1$ with trace field $K$ the vector space $K\cdot \Gamma$
naturally has the structure of a quaternion algebra $Q$ over $K$.
For each embedding $\sigma: K_1 \to \RR$ the quaternion
algebra $Q \times_{K_1,\sigma} \RR$ is isomorphic to the matrix ring 
$M_2(\RR)$ or to the Hamiltonians. In the first case the place 
$\sigma$ is said to be unramified, it the second case it is called
ramified.
\newline 
The analysis of the VHS of a Teichm\"uller curve yields:
\par
\begin{Cor} \label{quatalg}
The quaternion algebra $Q$ associated with $\SL(X,w)$
is isomorphic to $M_2(K)$. In particular all infinited places
are unramified and $\LL$ is defined over $K$, not only over $K_1$.
\end{Cor}
\par
{\bf Proof:}  
The result does not depend on passing from $\SL(X,\omega)$ to 
the subgroup of finite index $\Gamma$.
The fibre of $\oplus_{i=1}^r \LL_i$ over $c$ is defined over $\QQ$, 
say it equals $F \otimes_\QQ L$.
An element $\gamma \in \Gamma$ acts on this fibre as
$\diag(\gamma, \gamma^{\sigma_2}, \ldots, \gamma^{\sigma_r})$,
where $\sigma_i \in \Gal(L/\QQ)/\Gal(L/K)$.
This action commutes with the action of $k \in K$ on $F \otimes_\QQ L$ as
$\diag(k,k,k^{\sigma_2},k^{\sigma_2}, \ldots, k^{\sigma_r},k^{\sigma_r})$. 
Hence $K$ indeed acts on $F$ and we obtain a map 
$\psi: \Gamma \to \SL_2(K)$. 
The action on $F \otimes_\QQ \RR$ 
corresponds to postcomposing $\psi$ with the $r$ maps 
$\SL_2(K) \to \SL_2(\RR)$ from the embeddings $K \hookrightarrow \RR$.
The first embedding corresponds to the action on the
subspace $\langle \Re \omega, \Im \omega \rangle$. Thus the
Fuchsian embedding $\Gamma \subset \SL_2(\RR)$ we started with
factors through $\SL_2(K)$. 
Hence $K \cdot \Gamma = K \cdot \SL_2(K) = M_2(K)$.
\hfill $\Box$
\par

\subsection{Modular embeddings}
The expression 'modular embedding' is used in the literature
for maps from $\HH^n$ to the Siegel half space $\HH_n$ 
equivariantly for some group actions (see \cite{vG87}) 
and for maps $\HH \to \HH^n$ with some equivariance conditions
(see \cite{CoWo90}, \cite{ScWo00}). We compare that
latter type with maps  induced from Teichm\"uller curves.
\par
A subgroup $\Delta$ of $\SL_2(\RR)$ acting on
$\HH^n$ is called {\em arithmetic} if there is a quaternion
algebra $Q$ with $n$ unramified infinite places and an order
${\cal O}$ in $Q$ such that $\Delta$ is commensurable to the norm unit
group $\{ M \in {\cal O}: M{\cal O} \subset {\cal O}, \det(M) = 1 \}$.   
A Fuchsian group $\Gamma$ is said to have a {\em modular
embedding}  if there
exists an arithmetic group $\Delta$ acting on $\HH^n$ for
an appropriate $n$, an inclusion $\varphi: \Gamma \hookrightarrow \Delta$
and a holomorphic embedding 
$$\phi= (\phi_1,\ldots,\phi_n):\HH \to \HH^n$$
such that $\phi_1 = \id $ and $\phi(\gamma z) = \varphi(\gamma)\phi(z)$.
\par
\begin{Cor} \label{modemb}
Let $\Gamma$ be a subgroup of finite index in
the affine group $SL(X,\omega)$ of a Teichm\"uller curve, 
chosen as in Section \ref{univfam}. Then $\Gamma$ admits
a modular embedding $\HH \to \HH^r$, where $r=[K:\QQ]$.
\end{Cor}
\par
{\bf Proof:}
With the choices for the symplectic basis adapted to real multiplication
made in the setup (\ref{RM}), the universal property of the moduli
space of abelian varieties with real multiplication gives a
holomorphic map $\phi: \HH \hookrightarrow \HH^r$. Its components
are given by 
$$ \phi_i: \tau \mapsto \left(\int_{b_i} \omega_i(\tau) \right) / 
\left(\int_{a_i} \omega_i(\tau) \right), $$
where $\omega_i(\tau)$ are eigenforms for real multiplications
on the fibre of $A\to C$ pulled back to $\HH$.
The rest is parallel to the case $r=2$ (see \cite{McM1} Thm~10.1):
\newline
$\phi_1$ is an isometry (\cite{McM1} Thm.~4.1), hence a M\"obius
transformation. Replacing the generating flat surface $(X,\omega)$
by some $M\cdot(X,\omega)$ we may suppose it is the identity.
Consider a decomposition $A \hookrightarrow 
A_1 \times_C A_2$ up to isogeny. The endomorphisms of 
$A_1/C$ are an order ${\mathfrak o} \subset K$.
The moduli space of abelian varieties with real multiplication
by  ${\mathfrak o}$ is $\HH^r/\SL_2({\mathfrak o})$,
see Section \ref{RM}, and $\SL_2({\mathfrak o})$ is
an arithmetic subgroup of $\SL_2(\RR)$. Hence $\varphi$
descends to the moduli map for $A_1/C$ given by
$C \to \HH^r/\SL_2({\mathfrak o})$. This
implies that $\phi$ is equivariant with respect
to an inclusion $\varphi: \Gamma \to \Delta:=\SL_2({\mathfrak o})$,
which is of course the same as the map $\psi$ in
the proof of the previous 
corollary.
\hfill $\Box$
\par
Note that in the above proof the decomposition 
$A \hookrightarrow A_1 \times_C A_2$ up to isogeny is by no
means unique. Choosing a different one replaces the group
$\Delta$ by a commensurable arithmetic subgroup in $\SL_2(K)$.
\par
Corollary \ref{modemb} admits a converse, but depending on being able
to detect the image of $\mrC_g$ inside $\mrAV_g$ (Schottky problem).
\par
\begin{Thm} \label{Modemb}
Let $\Gamma$ be a Fuchsian group with totally real trace field
$K$ such that $g= [K:\QQ]$ and such that $\Gamma$ admits a modular embedding 
$\phi: \HH \to \HH^g$. If the composition of $\phi$ with the map to
the moduli space of abelian varieties, i.e.\ the map 
$$ \HH \to \HH^g \to \HH_g \to A_g$$
factors through the moduli space of curves $M_g$, then 
$\HH/\Gamma$ is a finite unramified cover of a 
Teichm\"uller curve in $\mrC_g$. 
\end{Thm}
\par
{\bf Proof:} If we lift the factorization through $\mrC_g$ 
to the universal covers we obtain maps
$$\HH \to T_g \to \HH_g \to \HH,$$ 
where the last map is $\HH_g \ni Z \mapsto z_{11}$. By definition
of a modular embedding the composition
of these maps is the identity, in particular an isometry
with respect to the Kobayashi (or equivalently Poincar\'e) metric.
Since a composition of maps is a Kobayahsi isometry if and only if
each single maps is, we conlcude that $\HH \to T_g$
is a geodesic. Since $\HH \to \HH_g$ is equivariant with respect to
$\Gamma \to \Delta$, the stabilizer of $\HH$ is a lattice and
$\HH \to T_g$ a Teichm\"uller curve.
\hfill $\Box$
\par
\subsection{A converse}
We now show that Teichm\"uller curves are characterized by a splitting of
the local system as in Prop.~\ref{splitasLvar}.
\par
\begin{Thm} \label{conversekrit}
Let $\pi: {\cal X} \to C = \HH/\Gamma$ be a family of curves of genus $g$
and suppose that the local system of the family $R^1 \pi_* \RR_{\cal X} $ has a
direct summand $\LL_\RR$ of rank $2$ given by the
Fuchsian embedding $i: \Gamma \hookrightarrow \SL_2(\RR)$ (up
to $\SL_2(\RR)$-conjugation).
\newline
Then $C \to \mrC_g$ is a finite covering of a Teichm\"uller curve.
\end{Thm}
\par
{\bf Proof:} By semisimplicity of the VHS (Thm.~\ref{Del87prop}) we have 
a decomposition of the local system 
$$ R^1 f_* \CC_{\cal X} =: \VV_\CC = \bigoplus_{i=1}^n (\LL_i \otimes W_i)$$
with $\LL_i$ irreducible local systems and $W_i$ vector spaces.
Since $\LL_\RR$ is defined via the Fuchsian embedding of $\Gamma$ it
is irreducible and we may suppose $(\LL_\RR) \otimes_\RR \CC \cong \LL_1$.
We do a priori not know that the $W_i$ are onedimensional.
But for a one-dimensional subspace $U \subset W_i$ of 
fixed bidegree Thm.\ \ref{Del87prop} shows that 
$\LL_1 \otimes U$ is a sub-$\CC$-VHS of $\VV_\CC$.
\newline
We consider the classifying map for $\LL_1 \otimes U$. In more concrete
terms, let $\omega_1(\tau)$ be a section over $\HH$ of the pullback of
the bundle $E^{1,0}(\LL_1 \otimes U) \subset  
E^{1,0}(\VV_C) = f_* \omega_{{\cal X}/C}$ to $\HH$ and complete it
to a basis $\omega_i(\tau)$, $i=1,\ldots,g$, of sections of the
pullback of $f_* \omega_{{\cal X}/C}$.
Choose a symplectic basis of sections $\{a_i,b_i\}$
of the pullback of $R^1 f_*  \CC_{\cal X}$ to $\HH$ such that
$\int_{a_i} \omega_j(\tau) = \delta_{ij}$ and 
$(\int_{b_i} \omega_j(\tau))_{i,j} \in \HH_g$. Then the classifying
map is given by the composition
$$ \phi: \HH \to T_g \to \HH_g \to \HH,$$
where the last map is $\HH_g \ni Z \mapsto z_{11}$.
\par 
The hypthesis states that $\phi$ is equivariant with respect 
to $\Gamma$ acting on domain and range via its (Fuchsian) 
embeddings to $\SL_2(\RR)$.
Hence $\phi$ descends to a holomorphic endomorphism $\ol{\phi}$
of $\HH/\Gamma$. Replacing $\Gamma$ by a subgroup of finite
index if necessary, we may suppose
$g(\HH/\Gamma) \geq 2 $. Since $\ol{\phi}$ cannot be constant 
it has to be an isomorphism. Thus $\phi$ is an isometry for
the Kobayashi metric. As in the previous theorem
this implies that $\HH \to T_g$ is an isometry. 
\hfill $\Box$
\par

\section{Higgs fields} \label{Higgssection}

Consider a weight one $\QQ$-VHS $(\VV, E^{1,0} \subset E=(\VV \otimes_\QQ
{\cal O}_C)_{\rm ext})$  
on a curve $C$ with unipotent monodromies. We have in mind that it comes from
a family of curves $f: {\cal X} \to C$. Recall that we 
denote the smooth compactification of $C$ by $\ol{C}$ and
we let $S = \ol{C}-C$. Let ${\cal X}_S:=\ol{f}^{-1}(S)$
the non-smooth fibres of the semistable family.
\newline
The connection $\nabla$ composed with the
inclusion and projection gives a map
$$ \Theta^{1,0}: E^{1,0} \to E \to E \otimes \Omega^1_{\ol{C}}(\log S) \to
(E/E^{1,0}) \otimes \Omega^1_{\ol{C}}(\log S) $$
that is ${\cal O}_{\ol{C}}$-linear. If we extend $\Theta^{1,0}$ by
zero mappings to the whole associated
graded sheaf $F:={\rm gr}(E) = E^{1,0} \oplus E^{0,1}$
of $E$ we obtain a  {\em Higgs bundle $(F,\Theta)$}. By definition this
is a vector bundle on $\ol{C}$ with a holomorphic map 
$\Theta: F \to F \otimes \Omega^1_{\ol{C}}(\log S)$, the additional
condition $\Theta \wedge \Theta = 0 $ being void if the base is a curve. 
\par
Sub-Higgs bundles of a Higgs bundles are subbundles, that are
stabilized by $\Theta$. 
Simpson's correspondence (\cite{Si90}) allows us to switch back and 
forth between sub-Higgs bundles of $F$ and sub-local systems
of $\VV$. 
\newline
The reason for working with Higgs bundles rather than
sublocal systems is twofold: In Lemma \ref{globalsec} it will be
useful to argue with degrees and the Higgs field has the simple
algebraic description as the edge morphism
$$ E^{1,0} = f_*\omega_{\ol{{\cal X}}/\ol{C}} 
\to R^1 f_*{\cal O}_{\cal X}
\otimes \Omega^1_{\ol{C}}(\log S) = E^{0,1} \otimes 
\Omega^1_{\ol{C}}(\log S)$$
of the tautological sequence
$$ 0 \to f^* \Omega^1_{\ol{C}}(\log S) \to \Omega^1_{\ol{{\cal X}}}(\log 
{\cal X}_S) \to \Omega^1_{\ol{{\cal X}}/\ol{C}}(\log {\cal X}_S) \to 0. $$
\par
Now for a Teichm\"uller curve $f: {\cal X} \to C$ we translate 
the decomposition in Prop.~\ref{splitasLvar} of $\VV_\CC := 
R^1 f_* \CC_{\cal X}$ in this language:
\newline
As in the general setup let $E$ the vector bundle associated 
with $\VV_\CC$ and $(F,\Theta)$ the corresponding Higgs bundle.
It decomposes into the sum of sub-Higgs bundles 
$$({\cal L}_i \oplus {\cal L}_i^{\vee}, \tau_i^{1,0}: 
{\cal L}_i \to {\cal L}_i^{\vee}
\otimes \Omega^1_{\ol{C}}(\log S)).$$
\newline
stemming from the $\LL_i$, $i=1,\ldots,r$ and the sub-Higgs 
bundle from $\MM$, which plays no role in the sequel. 
\par
\begin{Lemma} \label{HiggsofLi}
The Higgs field of all of the $\LL_i$ is
generically maximal i.e.\  $\tau^{1,0}_i$ is injective 
and precisely the Higgs field of $\LL_1$ is maximal, i.e.
$\tau_1^{1,0}$ is an isomorphism.
\end{Lemma}
\par
{\bf Proof:} By Prop.~4.10 in \cite{Ko85}, 
one can decompose a Higgs bundle
$F$ as a sum of Higgs bundles $F_a \oplus N$ with $N$ the maximal flat
subbundle, i.e.\ $\Theta^{1,0}|_N = 0$. 
For a weight one pVHS  the curvature of the Hodge bundle 
$E$ with respect to the Hodge metric coming from $H(\cdot,\cdot)$
can be expressed in terms of $\Theta^{1,0}$
(see \cite{Sd73} Lemma 7.18, compare also \cite{ViZu04} Section $1$).
Hence if $N$ is flat, the curvature vanishes and the bundle
$N$ comes from a unitary representation. 
\newline
We apply this to the Higgs bundles $({\cal L}_i \oplus {\cal L}_i^{\vee}, 
\tau_i^{1,0})$. 
The corresponding representation is Galois conjugate
to the Fuchsian representation of $\Gamma$, which contains 
non-trivial parabolic elements since $C$ is not compact. 
Thus the representation it is not unitary.
Hence $\tau^{1,0}$ is non-zero, proving the first statement.
\par
Lemma~2.1 in \cite{ViZu04} implies maximality for $\LL_1$.
If one of the $\LL_i$ for $i\neq 1$ was maximal, too, its Higgs
field would give an isomorphism ${\cal L}^{\otimes 2}_i \cong
\Omega^1_{\ol{C}}(\log S)$. Hence $\LL_i$ would become isomorphic
to $\LL_1$ after replacing $C$ by a finite cover.
This contradicts Lemma \ref{W1dim}, which is independent of passing
to a finite unramified cover.
\hfill $\Box$
\par
The converse implication of Lemma~2.1 in \cite{ViZu04} gives:
\par
\begin{Cor}
The groups $\Gamma^\sigma$ for $\id \neq \sigma \in \Gal(K/\QQ)$
are not Fuchsian.
\end{Cor}
\par
We let $d_i = \deg({\cal L}_i)$.
If we replace $C$ by an unramified covering of degree $n$ the
degrees $d_i$ will be multiplied by $n$. Hence the projectivized
$r$-tuple $(d_1:\ldots: d_r)$ seems to be an interesting invariant of 
the Teichm\"uller curve. Note that Lemma \ref{HiggsofLi} says
$d_1 > d_i$ for all $i$.
\subsection{Global sections}
Recall that the local system of a Teichm\"uller curve has a
decomposition $\VV_\QQ = \WW_\QQ \oplus \MM_\QQ$, and over $L$
we have the further decomposition
$\WW_L = \LL_1 \oplus \cdots \oplus \LL_r$. 
\par
\begin{Lemma} \label{globalsec11} 
The local system 
$\EEnd(\WW_L \otimes_L \CC)$ has no global sections of bidegree $(-1,1)$.
\end{Lemma}
\par
{\bf Proof:} As above
$$ \EEnd(\WW_\CC) = \bigoplus_{i=1}^r \EEnd((\LL_i)_\CC) \oplus 
\bigoplus_{1 \leq i,j \leq r, \, i\neq j} \HHom((\LL_i)_\CC,(\LL_j)_\CC). $$
Since the $\LL_i$ are irreducible and pairwise non-isomorphic, the
local system
$\HHom(\LL_i,\LL_j)$ has for $i \neq j$ no global sections
at all.
\newline
We show that $\EEnd(\LL_i)$ has only the global sections
$\CC \cdot \id$, which are in bidegree $(0,0)$. Recall that
$C$ is not compact, i.e.\ the group
$\Gamma$ contains a parabolic element $\gamma$. Its action
on a fibre of $\EEnd(\LL_i)$ decomposes into $2$ Jordan blocks,
one of size $1$ containing the identity and one of size $3$. Since
global sections are flat, they are contained in a unitary subbundle,
on which a unipotent element has to act trivially. Hence they
would have to correspond to the onedimensional eigenspace
of the Jordan block. But this contradicts the fact the flat
sections form a direct summand (again Prop.\ 4.10 in \cite{Ko85}).
\hfill $\Box$
\par
The following Lemma
will be used to determine Mumford-Tate groups in Section \ref{noconstant}.
The facts on Schur functors $\Schur_\lambda(\cdot)$ for a partition 
$\lambda$ of $\{1,\ldots,m\}$ used below 
can be found in \cite{FuHa91} Chapter $6$. 
\par
\begin{Lemma} \label{globalsec}
For all $(m,m')\in \NN^2$ the global sections of 
$\WW_\QQ^{\otimes m} \otimes (\WW^\vee_\QQ)^{\otimes m'}$ 
of bidegree $(0,0)$ are generated by tensor products
of global sections of the (trivial) bundles $\det(\LL_i)$.
\end{Lemma}
\par
{\bf Proof:} We tensor $\WW_\QQ$ and its summands by $\CC$, but
we omit this from notation. Renumber if necessary the indices such that
the $d_i$ are non-increasing. We have
$$\WW^{\otimes m} = \bigoplus_{(\lambda_1,\ldots,\lambda_k)} \; 
(\bigotimes_{i=1}^k 
\Schur_{\lambda_i}(\LL_i))$$
for partitions $\lambda_i$. The Schur functors for partitions with $3$ or more
rows are zero, because the $\LL_i$ are of rank
two, while we have for $a+b = \tilde{m}$, $a \geq b$
$$ \Schur_{\{a,b\}}(\LL_i) = \Sym^{\tilde{m}-2b}(\LL_i) 
\otimes \det(\LL_i)^a. $$
Since the $\LL_i$ have trivial $\det$-bundles and the $\LL_i$ are self-dual
it suffices  to show that 
$\widetilde{\LL} = \otimes_{i\in I}\, \Sym^{a_i}(\LL_i)$ 
has no global section of bidegree $(p,p)$ for
any non-empty subset $I \subset \{1,\ldots,r\}$ and any $a_i >0$ with
$\sum {a_i} = 2p$ and $1 \in I$. The last property maybe assumed
since we are interested in sections over $\QQ$ and the $\LL_i$
are all Galois conjugates. 
\newline
Suppose the contrary is the case.
Denote by
$(\widetilde{F},\theta)$ the Higgs bundle corresponding to $\widetilde{\LL}$.
The subsystem $\widetilde{\UU}$ of $\widetilde{\LL}$ generated by all
global sections is a polarizable sub-VHS, 
hence a direct summand by semisimplicity. By the hypothesis we
want to contradict there exists a non-trivial subbundle $\UU$ of 
$\widetilde{\UU}$ of type $(p,p)$. 
Since global sections are flat and the $(p,p)$-components
of flat section are again flat, the Higgs bundle associated
with $\UU$ is a direct summand of
$$\widetilde{F}^{p,p} = \bigoplus_{|\alpha| = p} \; F^\alpha,
\quad F^\alpha := \bigotimes_{i \in I} 
{\rm Sym}^{a_i} ({\cal L}_i \oplus {\cal L}_i^\vee)^{(\alpha_i,\,
a_i -\alpha_i)}, $$
where $\alpha = \{\alpha_i, i \in I \}$  and
$|\alpha| = \sum \alpha_i$.  
\newline
By Simpson's correspondence (\cite{Si90}) between local systems 
and Higgs bundles the sub-Higgs-bundle ${\cal U}$ is a direct 
summand and has degree $0$. This implies that 
$${\cal U} \hookrightarrow \bigoplus_{\deg \alpha =0} F^\alpha, $$
where $\deg \alpha = \sum_{i\in I} d_i (2\alpha_i - a_i)$ and $d_i = \deg
{\cal L}_i$.
\newline
Denote by $\Lambda$ the set of $\alpha^\lambda$ such that ${\cal U} \to
F^{\alpha^{\lambda}}$ is non-zero. Order $\Lambda$ lexicographically
starting with the $\alpha^\lambda$ where the first component
is maximal or if they are equal the second is, etc. 
Since $1 \in I$ we have $\alpha^\lambda_1 > 0$ for one $\alpha^\lambda$ 
hence by the ordering also $\alpha^1_1 >0$.
\newline
The Higgs field of $F^\alpha$ is the sum with appropriate signs 
$$ \Theta^{p,p}: F^\alpha \longrightarrow 
(\bigoplus_{\stackrel{|\beta|=1}{ \alpha - \beta > 0}} F^{\alpha-\beta}) 
\otimes \Omega^1_{\ol{C}}(\log S) $$
of the Higgs fields $\tau_i^{\alpha_i, a_i-\alpha_i}$ of the 
$\Sym^{a_i}(\LL_i)$. Since $d_1 >d_2$, 
since the $d_i$ are non-increasing  and since the 
$\alpha^\lambda \in \Lambda$ satisfy $\deg \alpha^\lambda = 0$
the subbundle $F^{\alpha -\beta} \otimes \Omega^1_{\ol{C}}(\log S)$
with $\alpha := \alpha^1$ and $\beta = (1,0,\ldots,0)$ 
only receives a Higgs field from $F^{\alpha^1}$ but from no
other $F^{\alpha^\lambda}$ for $\lambda \in \Lambda$. 
\newline
The Higgs field $\tau^{1,0}_1$ is a non-trivial by 
Lemma \ref{HiggsofLi} and hence the same is true for all the
$\tau_i^{\alpha_i, a_i-\alpha_i}$ . Since there
are no cancellations by the above argument, the composition
$$\xymatrix {
{\cal U} \ar[r] & \bigoplus_{\deg \alpha =0} F^\alpha \ar[r]^-{\Theta^{p,p}} 
& (\bigoplus_{\stackrel{|\beta|=1}{ \alpha - \beta > 0}} F^{\alpha-\beta}) 
\otimes \Omega^1_{\ol{C}}(\log S) 
\ar[r]^-{pr} 
& F^{\alpha^1 -(1,0,\ldots,0)} \otimes \Omega^1_{\ol{C}}(\log S)\\}$$
is non-trivial. This contradicts the flatness of $\UU$.
\hfill $\Box$
\par
\section{Shimura varieties} \label{noconstant}
Suppose we examine a Teichm\"uller curve where the local system
$\MM_\QQ$ of Prop.\ \ref{splitasLvar} vanishes, i.e.~such that $r=g$. 
This implies that $\VV$ does not split over $\QQ$. We will refer
to this condition by saying that the family of Jacobians
$A/C$ (or equivalently: its generic fibre) is simple. For dimension reasons 
the endomorphism ring of the generic fibre of $A$ cannot 
contain a bigger field than $K$. Nevertheless it might be possible
that a Teichm\"uller curve $C$ maps to a Shimura subvariety properly
contained in the locus of real multiplication. 
\newline
This Shimura subvariety could 
parametrize abelian varieties, whose endomorphism ring contains
a quaternion algebra but there are also Shimura subvarieties
that are not characterized by endormophisms of abelian
varieties, see \cite{Mu69}. 
\par
\subsubsection*{The Mumford-Tate group}

Let $B$ be an abelian variety with polarization $\Theta$ and
let $Q(\Theta)$ be the associated alternating form on $H^1(B,\QQ)$.
The {\em special Mumford-Tate group of $(B,Q$} is the largest 
algebraic $\QQ$-subgroup of ${\rm Sp}(H^1(B,\QQ), Q)$ that leaves 
the Hodge cycles of $B\times \cdots \times B$ invariant, i.e.\ 
all elements in $H^2p(B\times \cdots \times B,\QQ)^{p,p}$.
We denote it by $\Hg(B)$.
The {\em special Mumford-Tate group of
a polarized variation of Hodge structures $\WW_\QQ$} denoted by
$\Hg := \Hg(\WW_\QQ)$ is the largest $\QQ$ subgroup of 
${\rm Sp}(F, Q)$
fixing the Hodge cycles in $\WW^{\otimes m} \otimes
(\WW^{\vee})^{\otimes m'}$ of bidegree $(0,0)$, that remain Hodge under
parallel transform. Here $F$ is a fibre of $\WW_\QQ$ and
$Q$ the polarisation.
\newline
Mumford defines a {\em Shimura variety} (see \cite{Mu66}, \cite{Mu69})
as the coarse moduli space $M(\Hg)$ for the functor ${\cal M}(\Hg)$ 
representing isomorphism classes of polarized abelian varieties
with special Mumford-Tate group equal to a subgroup of $\Hg$.
Since the graphs of endomophisms of abelian varieties are Hodge 
classes, Shimura varieties of PEL-type (see \cite{Sh66}) fall within 
this category. 
\newline
The two definitions of Mumford-Tate groups are related by
the following general Lemma (see \cite{Sc96} Lemma 2.3, see also
\cite{Mo98} Section 1):
\par
\begin{Lemma}
If $\WW_\QQ = R^1 \pi_* \QQ$ for a family $g: B \to C$ of 
abelian varieties, then for a dense subset in $C$ the
special Mumford-Tate groups of the fibre coincide with $\Hg(\WW_\QQ)$. 
\end{Lemma}
\par
We specialize to a Teichm\"uller curve and to the family of
Jacobians $\pi: A \to C$:
\par
\begin{Thm} \label{Shimura}
Suppose that the family of abelian varieties $A/C$ over
a Teichm\"uller curve is simple or equivalently $[K:\QQ] = g$. 
\newline
Then the Shimura variety parametrizing abelian varieties with
real multiplication by $K$ is the smallest Shimura subvariety
of $\mrAV_g$ the Teichm\"uller curve $C$ maps to.
\end{Thm}
\par
{\bf Proof}: The locus of real multiplication is the
moduli space of the Mumford-Tate group $\Hg_K$ fixing the
Hodge cycles 
$$ \oplus_\sigma \; \sigma(a) \cdot \id_{\LL^ \sigma} \quad
{\rm for} \quad  a \in K.$$ 
By Lemma \ref{globalsec} the local system 
$\WW^{\otimes m} \otimes
(\WW^{\vee})^{\otimes m'}$ contains no  Hodge
cycles other than tensor powers and products of the above, 
hence $\Hg_K = \Hg(R^1 f_* \QQ).$ 
\hfill $\Box$
\par
\section{Rigidity, Galois action}
Teichm\"uller curves with $r=1$, i.e.\, with trace field $\QQ$ 
are obtained as unramified coverings of the once-punctured
torus by \cite{GuJu00} Thm.~5.5. 
They are called {\em square-tiled coverings} or {\em
origamis} because of this topological description. Since they
are Hurwitz spaces, they are known (\cite{Lo03}, \cite{Mo03})
to be defined over $\ol{\QQ}$. 
\newline
This also holds for all Teichm\"uller curves.
\par
Let $\WW$ be the local subsystem of $R^1 f_* \QQ_{\cal X}$
from Prop.\ \ref{splitasLvar}, which carries a polarized $\QQ$-VHS. 
The choice of a $\ZZ$-structure on $\WW$ defines a
family of abelian subvarieties $A_1/C \to A/C$. We could
take the families $A_1$ and $A_2$ of
Thm.\ \ref{Endallgemein} and dualize the isogeny to obtain
$$A_1 \hookrightarrow A_1 \times_C A_2 \to A.$$
In any case we do not claim any uniqueness of $A_1$.
We provide $A_1$ with the pullback of the polarisation
of $A$. It will be of type $\delta$, not necessarily
principal.  
\par
\begin{Thm} \label{nodefo}
For a Teichm\"uller curve $C$ and $A_1/C$ chosen as above
the canonical map $C \to A_{r,\delta}$ to the moduli
space of $\delta$-polarized abelian varieties of dimension $r$
admits no deformations. 
\newline
In particular the Teichm\"uller curve $C$ is defined over $\ol{\QQ}$.
\end{Thm}
\par 
{\bf Proof:}
By a theorem of Faltings (\cite{Fa83} Thm.~2) the tangent space to
the space of deformations
of $C \to A_{r,\delta}$ is a subspace of the 
global sections of $\End(\WW_\CC)$ of bidegree $(-1,1)$. The non-existence
of such sections was shown in Lemma \ref{globalsec11}. 
\newline
If the Teichm\"uller curve was not defined over $\QQ$ 
the transcendental parameters in the defining equations of $C$
would provide a non-trivial deformation of the family of 
Jacobians $A \to C$. \hfill $\Box$
\par
\begin{Rem} {\rm The proof of Thm.~\ref{Shimura} gives immediately
that the image of $C \to \mrAV_{r,\delta}$ lies in the locus of 
real multiplication but in no smaller Shimura variety. It would be 
interesting to know whether this is also the case for the map $C \to \mrAV_g$
if $[K:\QQ] < g$ and whether this map admits deformation.
}\end{Rem}
\par
\subsection{An algebraic description of Teichm\"uller curves}
The hypothesis of unipotent monodromies in the following theorem
can always be achieved by passing to a finite covering unramified
outside $S = \ol{C} - C$. It is there since 'maximal Higgs' 
does not make sense otherwise.
\par
\begin{Thm} \label{algchar}
Suppose that the Higgs bundle of family of curves 
$f: {\cal X} \to C = \HH/\Gamma$ with unipotent monodromies
has a rank two Higgs-subbundle with maximal Higgs field. 
Then $C \to \mrC_g$ is a finite covering of a Teichm\"uller curve.
\end{Thm}
\par
{\bf Proof:} Let $({\cal L} \oplus {\cal L}^{-1}, \tau^{1,0})$
be a maximal Higgs rank-$2$ subbundle.
By Simpsons correspondence (see \cite{Si90}
or the summary in \cite{ViZu04}) the local system $\VV_\CC$ 
has a rank two direct summand $\LL_\CC$ that carries a VHS, 
whose Higgs bundle is $({\cal L} \oplus {\cal L}^{-1}, \tau^{1,0})$.
\par
We claim that $\LL_\CC$ is defined over $\RR$: 
\newline
By the properties of a VHS the complex conjugate 
of $\LL_\CC$ has a Higgs bundle $({\cal L} \oplus {\cal L}^{-1}, 
\widetilde{\tau^{1,0}})$. Since the property maximal Higgs is 
equivalent to $2\deg {\cal L} = \deg \Omega^1_{\ol{C}}(\log S)$ we deduce that
$\widetilde{\tau^{1,0}}$ is also an isomorphism.
Moreover we must have $\LL_\CC \cong \ol{\LL_\CC}$ up to conjugation,
since otherwise the argument of \cite{McM1} Thm.~4.2 explained
in Remark \ref{McM42Rem} leads to a contradiction. In particular
the traces of $\LL_\CC$ are real.
\newline
Now consider $\LL_\CC$ as given by a representation 
$$\rho: \pi_1(C,c) \to \SL_2(\CC).$$ 
We cannot yet use that the image of $\rho$ equals $\Gamma$. 
But we know that $\rho$ is semisimple and irreducible even
if we pass to a subgroup of finite index in $\pi_1(C,c)$. 
From this one deduces that the image of $\rho$ contains
two non-commuting hyperbolic elements. If the eigenvalues
of one such element are real, we are done by \cite{Ta69} (compare
Lemma \ref{fromTake}).
\newline
So suppose all these eigenvalues are non-real.
If $\rho$ and $\ol{\rho}$ are given in the
eigenbasis of a hyperbolic element $\gamma_1$, a matrix $M$ that
conjugates $\rho$ into $\ol{\rho}$ has to be off diagonal.
If $\gamma_2$ is hyperbolic and does not commute with $\gamma_1$
they do not share an eigenvector. Since $M$ has to be
off diagonal, too, if we consider $\rho$ and $\ol{\rho}$ in the
eigenbasis for $\gamma_2$, we obtain a contradiction.
\par
Using the claim and by Lemma 2.1 in \cite{ViZu04} the
property 'maximal Higgs' implies that the action 
of $\pi_1(C)$ on a fibre of $\LL$ is just the action of 
$\Gamma \subset \SL_2(\RR)$. 
We can now use Thm.\ \ref{conversekrit} to conclude.
\hfill $\Box$
\par
This characterization may not be very useful to construct
Teichm\"uller curves but it has the advantage of being
completely algebraic. In \cite{Mo03} we noted that
the absolute Galois group $G_\QQ$ acts on the set of
all origamis (in fact faithfully in an appropriate sense). 
\par
\begin{Cor} \label{GQacts}
The absolute Galois group acts on the set
of Teichm\"uller curves with $r=g$.
\end{Cor}
{\bf Proof:} In case $r=g$ the map $C \to M_g$ is defined
over a number field, since $C \to A_g$ is by Thm.~\ref{nodefo}.
Hence it makes sense to let $G_\QQ$ act on the family 
$f: {\cal X} \to C$.
\newline
The construction of the Higgs bundle (see Section $3$) is
algebraic. So the Higgs bundle of the Galois conjugate curve 
will have as many Higgs subbundles of a given rank as the
original one. And the $G_\QQ$-action on families of
curves defined over number fields preserves the property of
the Higgs field to be an isomorphism on a subbundle.
\hfill $\Box$
\par
\begin{Rem} {\rm 
For the cases $1<r<g$ we cannot show with these methods
that the map $C \to M_g$ is defined over a number field,
although we expect this to be true in general. 
\newline
If it is true, then the above proof applies to show that
$G_\QQ$ acts on the set of all Teichm\"uller curves.
}\end{Rem}


Martin M{\"o}ller: Universit{\"a}t Essen, FB 6 (Mathematik) \newline 
45117 Essen, Germany \newline
e-mail: martin.moeller@uni-essen.de \newline

\end{document}